\documentclass[a4paper,10pt]{article}

\usepackage{amsmath,amssymb,amsthm}

\def\la{\lambda}
\def\g{\mathfrak{g}}
\def\b{\mathfrak{b}}
\def\g{\mathfrak{g}}
\def\h{\mathfrak{h}}
\def\n{\mathfrak{n}}
\def\l{\mathfrak{l}}
\def\N{\mathbb{N}}
\def\Z{\mathbb{Z}}
\def\C{\mathbb{C}}
\def\F{\mathbb{F}}
\def\qed{$\hfill \blacksquare$}

\newtheorem{theo}{Theorem}
\newtheorem{prop}{Proposition}
\newtheorem{lem}{Lemma}
\numberwithin{equation}{section}
\numberwithin{lem}{section}
\numberwithin{theo}{section}
\numberwithin{prop}{section}

\begin{document}
\title{Twisted traces of intertwiners for Kac-Moody algebras and
 classical dynamical r-matrices corresponding to generalized Belavin-Drinfeld 
triples}

\author{Pavel Etingof and Olivier Schiffmann}
\date{}
\maketitle

\section{Introduction}
In early eighties, Belavin and Drinfeld \cite{BD}
classified nonskewsymmetric classical r-matrices
for simple Lie algebras. It turned out that such r-matrices,  
up to isomorphism and twisting by elements from 
the exterior square of the Cartan subalgebra, are classified by 
rather unusual combinatorial objects which are now 
called Belavin-Drinfeld triples. By definition, a Belavin-Drinfeld triple 
for a simple Lie algebra $\g$ is 
a triple $(\Gamma_1,\Gamma_2,T)$, where $\Gamma_1,\Gamma_2$ are subsets 
of the Dynkin diagram $\Gamma$ of $\g$, and $T:\Gamma_1\to \Gamma_2$ is 
an isomorphism of Dynkin diagrams, which satisfies the nilpotency 
condition: if $\alpha\in \Gamma_1$ then there exists $k$ such that 
$T^{k-1}(\alpha)\in \Gamma_1$ but $T^k(\alpha)\notin \Gamma_1$. 
The r-matrix corresponding to such a triple is given by 
a certain explicit formula. This formula works not only for 
simple finite dimensional Lie algebras but in fact for any symmetrizable 
Kac-Moody algebra. 

In \cite{S}, the second author generalized the work of Belavin and Drinfeld 
and 
classified  classical
nonskewsymmetric \textit{dynamical} r-matrices for simple Lie algebras. 
It turns out that they have an even simpler classification: up 
to gauge transformations, they are classified by generalized 
Belavin Drinfeld triples, which are defined as the usual Belavin-Drinfeld 
triples but without any nilpotency condition.  
The dynamical r-matrix corresponding to such a triple is given by 
a certain explicit formula, which, as before, works not only for 
simple finite dimensional Lie algebras but in fact for any symmetrizable 
Kac-Moody algebra. 

This includes
some well known examples: if $T=id$ one gets Felder's dynamical r-matrix, and if $\g$ is
of type 
$\hat A_{n-1}$ (i.e. the Dynkin diagram is an n-gon) 
and $T$ is the rotation of the n-gon by an angle $2\pi k/n$ where $k$ is 
prime to $n$ then one gets Belavin's elliptic r-matrix. 

G.Felder \cite{Fel} associated 
to every classical dynamical r-matrix, 
a remarkable 
system of differential equations called the Knizhnik-Zamolodchikov-Bernard
(KZB) equations. For the Felder and the Belavin r-matrix 
these equations have a 
representation-theoretical interpretation. Namely, the 
KZB equations with the Felder r-matrix are satisfied by  
conformal blocks for the Wess-Zumino-Witten (WZW) model of conformal field theory on elliptic curves 
(or, equivalently, weighted traces of products of intertwining operators)
\cite{Ber},\cite{Fel}, and the 
KZB equations with 
the Belavin r-matrix are satisfied by conformal blocks for the 
WZW model twisted by the rotation 
of the Dynkin diagram \cite{Et}, \cite{Tak}. It is therefore natural to expect 
that a similar interpretation exists for all dynamical r-matrices.

In this paper, we propose such an interpretation. Namely, 
for any Kac-Moody algebra and a (nondegenerate) generalized BD triple we 
show that weighted 
traces of products of intertwining 
operators, "twisted" by $T$, satisfy 
the KZB equations with the corresponding dynamical r-matrix from \cite{S}. 

We consider two cases: operators with values 
in representations from category $\cal O$, and operators 
with values in finite dimensional representations for affine Lie algebras. 
In the first case we get  the KZB equations 
for trigonometric dynamical r-matrices 
corresponding to generalized Belavin-Drinfeld 
triples for Kac-Moody algebras, and 
in the second case we get 
the KZB equations for elliptic dynamical r-matrices which 
are intermediate between Felder's and Belavin's elliptic r-matrices. 

In conclusion we would like to point out some directions 
of future research. 

First of all, it turns out that classical dynamical
r-matrices mentioned above can be explicitly quantized. 
To obtain such a quantization has been an open problem, except 
for a few special cases, but we will present a complete solution of this 
problem in our next paper (joint with Travis Schedler). 
 
In another
 paper we plan to generalize the results of the present paper 
to quantum groups, following the ideas
of \cite{EV3}, where it was done for the case $\Gamma_1=\Gamma_2=\Gamma$ and $T=\mathrm{Id}$.
This will give quantum KZB equations, 
which are difference equations involving quantum dynamical R-matrices which are quantizations of the classical dynamical r-matrices
that appeared in this paper. In the case of the Belavin R-matrix, 
these equations are the elliptic qKZ equations for the 8-vertex 
model, which play an important role in statistical mechanics. 

\section{Classical dynamical r-matrices for Kac-Moody algebras}
\paragraph{}In \cite{S} the second author associated a solution of the classical dynamical
Yang-Baxter equation to every generalized Belavin-Drinfeld triple for simple Lie algebras.
This construction easily extends to any Kac-Moody algebra. We recall this construction
in this section. 
\paragraph{Preliminaries.} Let $A=(a_{ij})$ be a symmetrizable
generalized Cartan matrix of size $n$ and rank $l$. Let
$(\h,\Gamma,\check{\Gamma})$ be a realization of $A$, i.e $\h$ is a complex
vector space of dimension $2n-l$, $\Gamma=\{\alpha_1,\ldots
\alpha_n\}\subset \h^*$ and $\check{\Gamma}=\{h_1,\ldots h_n\}\subset \h$
are linearly independent sets and $\langle \alpha_j,h_i\rangle =a_{ij}$.
Let $\g=\n_- \oplus \h \oplus \n_+$ be the Kac-Moody algebra associated
to $A$, i.e $\g$ is generated by elements $e_i$, $f_i$, $i=1,\ldots n$ and
$\h$ with relations
$$[e_i,f_j]=\delta_{ij}h_i,\qquad [\h,\h]=0,\qquad [h,e_i]=\langle
\alpha_i,h\rangle e_i, \qquad [h,f_i]=-\langle \alpha_i,h\rangle f_i,$$
together with the Serre relations (see \cite{K}). Set $\b_\pm=\h \oplus \n_\pm$. Let
$\h'=\bigoplus
\C h_i$ and let $(\,,\,)$ be a nondegenerate invariant bilinear
form on $\g$. Recall that the restriction of $(\,,\,)$ to $\h$ is
nondegenerate (hence it defines a form on $\h^*$, which we also denote 
by $(,)$), and that the kernel of the restriction of $(\;,\;)$ to $\h'$
is equal to the center $\mathfrak{c}$ of $\g$. Let $\Delta=\Delta^+ \cup
\Delta^- \subset \h^*$ be the root system of $\g$ and let $\g_\alpha$
denote the weight subspace corresponding to a root $\alpha$. For each
$\alpha \in \Delta^+$ fix bases $e_{\alpha}^{(1)},\ldots
e_{\alpha}^{(dim\;\g_\alpha)}$ and $f_{\alpha}^{(1)},\ldots 
f_{\alpha}^{(dim\;\g_\alpha)}$ of $\g_\alpha$ and $\g_{-\alpha}$
respectively such that $(e_{\alpha}^{(i)},f_{\alpha}^{(j)})=\delta_{ij}$.
The vector $[e_\alpha^{(i)},f_\alpha^{(i)}]$ is independent of the choice of $i$ and we set $h_\alpha=[e_\alpha^{(i)},f_\alpha^{(i)}]$ for any
$i=1\ldots \mathrm{dim}\;\g_\alpha$. Let $\rho \in \h^*$ be an 
element satisfying $(\rho,h_i)=1$ for all $i$. Identifying $\h$ with $\h^*$ via the form
$(\,,\,)$, we can regard $\rho$ as an element of $\h$ ($\rho$ is
well-defined up to adding a central element).
\paragraph{Definition.} A \textit{generalized Belavin-Drinfeld triple} is a triple
$(\Gamma_1,\Gamma_2,T)$ where $\Gamma_1,\Gamma_2 \subset \Gamma$ and $T:
\Gamma_1 \stackrel{\sim}{\to} \Gamma_2$ is a inner product 
preserving bijection.
\paragraph{}Given a generalized Belavin-Drinfeld triple $(\Gamma_1,\Gamma_2,T)$, we
let $\Gamma_3$ be the subset of $\Gamma_1 \cap \Gamma_2$ 
consisting of roots which return to their original position after applying 
$T$ several times. It is clear that $(\Gamma_1\setminus\Gamma_3,\Gamma_2\setminus\Gamma_3,T)$ is an ordinary Belavin-Drinfeld triple 
(i.e. $T$ satisfies the nilpotency condition) and
$(\Gamma_3,\Gamma_3,T)$ is generalized Belavin-Drinfeld triple on which $T$ is an automorphism of the Dynkin subdiagram $\Gamma_3$. 
\paragraph{}Set $\l=\big(\sum_{\alpha \in \Gamma_1} \C(\alpha
-T(\alpha))\big)^\perp\subset \h$ and let $\h_0\subset \h$
be the orthogonal complement of $\l$ in $\h$ with respect to the 
inner product on $\h$. 

\paragraph{Definition.} We say that the generalized Belavin-Drinfeld triple is
\textit{nondegenerate} if the restriction of $(\,,\,)$ to $\l$ is nondegenerate. 

In the nondegenerate case, we have $\h=\l\oplus \h_0$. 

\paragraph{}Note that every generalized Belavin-Drinfeld triple on a Dynkin diagram of finite type or a connected Dynkin diagram of affine type
is nondegenerate. Indeed, in a finite type case, 
 $\l$ is a real subalgebra of $\h$, and 
the form on the real part of $\h$ is positive definite, 
so $\l$ is nondegenerate. A similar argument works in the affine case.  
Nondegeneracy  is also the case when $\Gamma_1=\Gamma_2$, 
because in this case $T$ gives rise to a finite order (hence semisimple) 
orthogonal automorphism of $\h$, and $\l=\h^{T}$. 

\paragraph{}Let $(\Gamma_1,\Gamma_2,T)$ be a nondegenerate generalized
Belavin-Drinfeld triple. It is convenient to choose an orthonormal basis
 $(x_j)_{j \in I}$ of $\h$ with respect to $(\;,\;)$ in such a way that 
$$\l=\bigoplus_{j \in I_1} \C x_j , \qquad \h_0=\bigoplus_{j \in I_2} \C x_j,$$
for suitable disjoint subsets $I_1,I_2 \subset I$ such that $I=I_1 \cup I_2$.

Let $\h_i$ be the subspaces of $\h$ spanned by $h_\alpha,\alpha\in \Gamma_i$. 

The following Lemma is straightforward but important for 
the considerations below. 

\begin{lem} There exists a unique Lie algebra homomorphism $B:\n_-\oplus \h_1
\to\n_-\oplus \h_2$ (resp.
$B^{-1}: \n_+\oplus \h_2 \to \n_+\oplus \h_1$) such that
$B(f_\alpha)=f_{T(\alpha)}$, $B(h_\alpha)=h_{T(\alpha)}$ if $\alpha \in \Gamma_1$,
$B(f_\alpha)=0$ if
$\alpha
\not\in \Gamma_1$ (resp.
$B^{-1}(e_\alpha)=e_{T^{-1}(\alpha)}$, $B^{-1}(h_\alpha)=h_{T^{-1}(\alpha)}$ if
$\alpha \in \Gamma_2$, $B^{-1}(e_\alpha)=0$ if
$\alpha
\not\in \Gamma_2$).\end{lem}

We will extend the homomorphisms $B,B^{-1}$ of Lie algebras to the 
corresponding homomorphisms of their universal enveloping algebras. 

\paragraph{}Let $\langle \Gamma_i \rangle$,
$i\in \{1,2,3\}$ be the set of roots
$\alpha
\in \Delta^+$ which are linear combinations of simple roots from $\Gamma_i$. Let
$\g_{\Gamma_i}$ be the subalgebra of $\g$ generated by $\g_\alpha,\g_{-\alpha}$,
$\alpha \in \Gamma_i$. The map $B$ restricts to an automotphism of $\g_{\Gamma_3}$. For
each
$\alpha
\in
\langle\Gamma_3\rangle$, let $N_\alpha$ be the order of the action of $B$ on $\alpha$, i.e
$B^{N_\alpha}\alpha=\alpha$ but $B^r\alpha \neq \alpha$ for \linebreak $1 \leq r <
N_\alpha$.  Finally, it will be convenient to assume
that $f_\alpha^{(t)}$ and $e_\alpha^{(t)}$ are eigenvectors for
$B^{N_\alpha}$, and set $B^{N_\alpha} f_\alpha^{(t)}=\theta_\alpha^{(t)}
f_\alpha^{(t)}$. Note that
$B^{N_\alpha}e_\alpha^{(t)}=\theta_\alpha^{(t)-1}e_\alpha^{(t)}$.

Finally, let us define the Cayley transform of $T$ on $\h_0$. 
To do this, we need the following straightforward Lemma.

\begin{lem} For any $x\in \h_0$, there exists a unique $y\in \h_0$ such that 
for all $\alpha\in \Gamma_1$ one has $(\alpha-T\alpha,y)=(\alpha+T\alpha,x)$. 
\end{lem}

It is clear that $y$ depends linearly on $x$. We will write 
$y=C_Tx$. It is easy to check that the operator 
$C_T:\h_0\to \h_0$ is skewsymmetric. It is called the Cayley 
transform of $T$.  

\begin{prop}[\cite{S}]\label{P:S1} The function $r_T: \l^* \to (\g \otimes \g)^\l$
\begin{equation*}
\begin{split}
r_T(\lambda)=& -\frac{1}{2}\sum_{j \in I} x_j \otimes x_j
+\frac{1}{2}\sum_{i \in I_2}C_Tx_i\otimes x_i - \sum_{\alpha,t}
f_\alpha^{(t)} \otimes e_\alpha^{(t)}
\\
&+\sum_{\alpha,t} \sum_{l=1}^\infty
e^{-l(\alpha,\lambda)} e_\alpha^{(t)} \wedge B^lf_\alpha^{(t)}
\end{split}
\end{equation*}
is a solution of the classical dynamical Yang-Baxter equation
\begin{equation*}
\begin{split}
\sum_i &\left(x_i^{(1)} \otimes \frac{\partial}{\partial
x_i}r^{23}(\lambda)-x_i^{(2)}
\frac{\partial}{\partial x_i}r^{13}(\lambda)+x_i^{(3)}\frac{\partial}{\partial
x_i}r^{12}(\lambda)\right)\\
&+[r^{12}(\lambda),r^{13}(\lambda)]+[r^{13}(\lambda),r^{23}(\lambda)]+[r^{12}(\lambda),r^{23}(\lambda)]=0.
\end{split}
\end{equation*}
\end{prop}
\paragraph{Remarks.} i) In the expression for
$r_T(\lambda)$, the sum
$\sum_{l=1}^\infty e^{-l(\alpha,\lambda)}e_\alpha^{(t)} \wedge B^lf_\alpha^{(t)}$ is finite
if $\alpha
\not\in \langle\Gamma_3\rangle$ and is 
an infinite series convergent to 
a rational function of $e^{(\alpha,\lambda)}$ if
$\alpha \in \langle \Gamma_3 \rangle$.

Consider the special case when $\Gamma_1=\Gamma_2=\Gamma$, and $T$ is 
an automorphism of the Dynkin diagram. Let $N$ be the order 
of $B$ (i.e. $N$ is the smallest number divisible by all $N_\alpha$ 
such that  
$(\theta_\alpha^{(t)})^{N/N_\alpha}=1$).
In this case, the formula for $r_T(\la)$ can be written in the form
\begin{equation*}
\begin{split}
r_T(\lambda)=& -\frac{1}{2}\sum_{j \in I} x_j \otimes x_j+\frac{1}{2}\sum_{i \in
I_2}C_Tx_i\otimes x_i
 - \sum_{\alpha,t}
f_\alpha^{(t)} \otimes e_\alpha^{(t)}\\
&+\sum_{\alpha,t} \sum_{l=1}^{N}
\frac{e^{-l(\alpha,\lambda)}}{1-e^{-N(\alpha,\la)}} 
e_\alpha^{(t)} \wedge B^lf_\alpha^{(t)}. 
\end{split}
\end{equation*}

\paragraph{Examples.} 
i) When
$T=1$ one obtains Felder's trigonometric dynamical r-matrix
$$r(\lambda)=-\frac{\Omega}{2}+\sum_{\alpha>0,t}\frac{1}{2}\mathrm{cotanh}(\frac{1}{2}(\alpha,\lambda))
e_\alpha^{(t)} \wedge
f_\alpha^{(t)},$$ where $\Omega=\sum_{\alpha,t} (e_\alpha^{(t)}
\otimes f_\alpha^{(t)} + f_\alpha^{(t)} \otimes e_\alpha^{(t)}) + \sum_i x_i \otimes x_i$ is
the Casimir element (see \cite{EV1}, Section 3).

ii) Consider $\g=sl(3)$ and let $T$ be the automorphism exchanging
the two simple roots. Then $\l=\C \rho$, so we can regard the element
$\la\in \l^*$ as a scalar. In this case, 
the dynamical r-matrix $r_T(\lambda)$ is 
\begin{equation*}
\begin{split}
r_T(\lambda)=&
-\frac{1}{2}\sum_{i \in I} x_i \otimes x_i-
\frac{e^{-2\lambda}}{1-e^{-2\lambda}}(e_{\alpha_1}\otimes
f_{\alpha_1}+e_{\alpha_2}\otimes f_{\alpha_2}) 
-\frac{e^{-2\lambda}}{1+e^{-2\lambda}}e_{\alpha_1+\alpha_2}\otimes
f_{\alpha_1+\alpha_2}\\
 &-\frac{1}{1-e^{-2\lambda}}(f_{\alpha_1} \otimes e_{\alpha_1}+f_{\alpha_2} \otimes
e_{\alpha_2}) -\frac{1}{1+e^{-2\lambda}}f_{\alpha_1+\alpha_2} \otimes
e_{\alpha_1+\alpha_2}\\
 &+\frac{e^{-\lambda}}{1-e^{-2\lambda}}(e_{\alpha_1} \wedge
f_{\alpha_2} +e_{\alpha_2}\wedge f_{\alpha_1}).
\end{split}
\end{equation*}

Note that in this example the eigenvalue $\theta_{\alpha_1+\alpha_2}^{(1)}$
of $B$ is $-1$, which is the reason for the appearance of denominators 
$1+e^{-2\la}$.

\paragraph{}In the next section we give the representation-theoretic interpretation of the KZB equations associated to the dynamical 
r-matrix $r_T(\lambda)$.
\section{Traces of intertwining operators}
\paragraph{3.1. Traces.}

For any $\h$-diagonalizable $\g$-module $V$ let $V[\lambda]$ denote the subspace of $V$ of weight $\lambda \in \h^*$. Let $M_\lambda$ be the Verma module of highest weight $\lambda$ and let $v_\lambda \in M_\lambda[\lambda]$ be a highest weight vector. Let $M_\lambda^*$ be the graded dual Verma module: $M_{\lambda}^*=\bigoplus_{\mu}M[\mu]^*$ as a vector space and the Lie algebra $\g$ acts by
$$x.a(u)=-a(x.u) \qquad \forall\;\;x \in \g,\;u \in M_\lambda,\;a \in M_{\lambda}^*.$$
Let $v_\lambda^* \in M_{\lambda}^*[-\lambda]$ be the lowest weight vector satisfying $\langle v_\lambda^*,v_\lambda \rangle=1$. Note that $M_\lambda$ and $M_\lambda^*$ are irreducible for generic values of $\lambda$.

\paragraph{}Recall the definition of the quadratic Casimir operator
$C$: 
$$C=2\rho + \sum_j x_j^2 + 2\sum_{\alpha \in \Delta^+}\sum_i f_{\alpha}^{(i)}e_{\alpha}^{(i)}.$$
The operator $C$ acts on $M_\lambda$ by multiplication by $\Delta_\lambda:=(\lambda,\lambda+2\rho)$. 

\paragraph{Intertwining operators.} Let $\lambda,\mu \in \h^*$ and let $V$ be a $\g$-module from the category $\mathcal{O}$. We will consider compositions of intertwining operators of the form
$$\Phi:\; M_\lambda \to M_\mu \otimes V.$$
The following lemma is well known (a proof can be found e.g in \linebreak \cite{EV2},
\cite{EFK1}) :

\begin{lem} Suppose that $M_\mu^*$ is irreducible. Then the map
$$\mathrm{Hom}_\g(M_\lambda, M_\mu \otimes V) \to V[\lambda-\mu],\; \Phi \mapsto \langle v_\mu^*, \Phi v_\lambda\rangle$$
is an isomorphism. \end{lem}

\paragraph{}Given $v \in V[\lambda-\mu]$, and when $M^*_\mu$ is irreducible, we will denote by $\Phi^v_\lambda$ the unique intertwiner $M_\lambda \to M_\mu \otimes V$ satisfying $\langle v_\mu^*, \Phi v_\lambda
\rangle=v$. 
\paragraph{}For any $\lambda \in \h^*$ and any $\h$-semisimple $\g$-module $V$ we denote by $\lambda \in \mathrm{End}\;(V)$ the operator satisfying
$\lambda_{|V[\nu]}=(\lambda,\nu)$.
\paragraph{} Let $\mu,\mu' \in \h^*$.
Consider the linear operator $B: M_{\mu'} \to M_\mu$ defined by $B(xv_{\mu'})=B(x)
v_\mu$ for any $x \in U\n_-$. The following lemma is straightforward.
\begin{lem}\label{Blem}Let ${\mu},\mu'$ be such that 
$(\mu',\alpha)=(\mu,T\alpha)$ for all 
$\alpha\in \Gamma_1$. Then for every $x \in U(\n_-\oplus \h_1)$ 
we have $Bx=B(x)B$ and for
every $x \in U(\n_+\oplus\h_2)$ we have $xB=BB^{-1}(x)$. \end{lem}

\paragraph{}Let
$V_1,V_2,\ldots V_r$ be $\g$-modules from the category $\mathcal{O}$,
$v_1 \in V_1, \ldots, v_r \in V_r$ homogeneous vectors of weights
$\nu_1,\ldots \nu_r$ respectively. Set $\nu=\sum \nu_i$ and consider
(for generic $\mu$) the composition
\begin{equation}\label{EE:1}
\Phi^{v_1}_{\mu-\nu_2-\ldots -\nu_{r}} \ldots \Phi^{v_r}_{\mu} B e^\lambda: M_{\mu'} \to
 M_{\mu-\nu} \otimes V_1 \otimes \ldots \otimes V_r 
\end{equation}
where $\lambda \in \l^*,\;\mu' \in \h^*$, $(\mu',\alpha)=(\mu,T\alpha)$
for $\alpha\in \Gamma_1$
(we regard $\l^*$ as a subspace of $\h^*$ using the inner product).
\paragraph{}If
\begin{equation}\label{E:1}
\mu'=\mu-\nu
\end{equation}
we can define
\begin{equation*}
F^{v_1,\ldots,v_r}(\lambda,\mu)=\mathrm{Tr}\;(\Phi^{v_1}_{\mu-\nu_2-\ldots
-\nu_{r}} \ldots \Phi^{v_r}_{\mu} B e^\lambda): (\h^*)^B \to V_1
\otimes \ldots \otimes V_r. 
\end{equation*}
If $\nu\in\l^\perp$, 
the space of solutions of (\ref{E:1}) in $\mu$ is an
$\l^*$-principal homogeneous space. Note that it
follows from the Kac-Kazhdan conditions (see \cite{KK}) and from the fact that $\rho
\in \l$ that for any fixed $v_1,\ldots,v_r$ the composition (\ref{EE:1}) is defined for
generic values of $\mu$ in any $\l^*$-principal homogeneous space.
In particular, $F^{v_1,\ldots ,v_r}(\lambda,\mu)$
is a formal series in $\lambda$ whose coefficients are 
trigonometric functions of $\mu$ with values in
the space $(V_1 \otimes \ldots  \otimes V_r)^\l$. 

\paragraph{}Set
$$\delta_B(\lambda)=\big(\mathrm{Tr}_{|M_{-\rho}}(
Be^\lambda)\big)^{-1}$$
(a twisted
version of the Weyl denominator).
\begin{lem} We have $\delta_B(\lambda)=e^{(\rho,\lambda)}\prod_{\overline{\alpha} \in
\langle \Gamma_3 \rangle/B}\prod_{t} (1-\theta_\alpha
^{(t)}e^{-N_\alpha(\alpha,\lambda)})$. \end{lem}
\noindent
\textit{Proof.} Identify $U\n^-$ with $M_{-\rho}$ by $u \mapsto uv_{-\rho}$. Let
$\C=U_0 \subset U_1 \subset \ldots$ denote the canonical filtration of $U\n^-$. Since
$(U_n)$ is stable under $Be^\lambda$, we can replace $U\n^-$ by its graded when computing
$\mathrm{Tr}_{|U\n^-}(Be^\lambda)$. By the PBW theorem, $\mathrm{Gr}\;U\n^-=S \n^-$. By
definition $B$ acts nilpotently on $f_\alpha^{(t)}$ unless $\alpha \in \langle
\Gamma_3\rangle$. Hence
$\mathrm{Tr}_{|S\n^+}(Be^\lambda)=\mathrm{Tr}_{|S\n^-_{\Gamma_3}}(Be^\lambda)$
where $\n^-_{\Gamma_3}=\n^-\cap \g_{\Gamma_3}$. The Lemma now follows from
$$S\n^-_{\Gamma_3}=\bigotimes_{\overline{\alpha} \in
\langle \Gamma_3 \rangle/B,t} S(\C f_\alpha^{(t)} \oplus \ldots \oplus
\C B^{N_\alpha-1}f^{(t)}_\alpha)$$
and $$\mathrm{Tr}_{|S(\C f_\alpha^{(t)} \oplus \ldots \oplus
\C B^{N_\alpha-1}f^{(t)}_\alpha)}(Be^\lambda)=\frac{1}{1-\theta_\alpha
^{(t)}e^{-N_\alpha(\alpha,\lambda)}}.$$
\qed

 We put $\F^{v_1,\ldots
,v_r}(\lambda,\mu)=\delta_B(\lambda)F^{v_1,\ldots ,v_r}(\lambda,\mu)$.
This function is the main object of this paper. 

\paragraph{3.2. The KZB equations for traces.}

The following theorem is one of our main results.
 
\begin{theo} The function $\F^{v_1,\ldots ,v_r}(\lambda,\mu)$ satisfies the following
system of differential equations, for $i=1,\ldots r$:
\begin{equation}\label{E:prop1}
\begin{split}
\bigg(\sum_{j \in I_1} x_{j|V_i} \frac{\partial}{\partial
x_j}+&\sum_{j>i}{r}_T(\lambda)_{|V_i \otimes V_j} -
\sum_{j<i}{r}_T(\lambda)_{|V_j \otimes V_i} \bigg) \F^{v_1,\ldots
,v_r}(\lambda,\mu)\\
&=\frac{1}{2}\left(\Delta_{\mu-\nu_{i+1}-\ldots-\nu_r}-\Delta_{\mu-\nu_{i}-\ldots
-\nu_{r}}\right)\F^{v_1,\ldots ,v_r}(\lambda,\mu)
\end{split}
\end{equation}
\end{theo}
\noindent
These equations are called the KZB equations (see \cite{Fel}). 

\noindent
\textit{Proof.} For simplicity of notation we will write
$\mu_i=\mu-\nu_r-\ldots -\nu_{i+1}$. We compute the function
$$A(\lambda,\mu)=\mathrm{Tr}(\Phi^{v_1}_{\mu_{1}}\ldots
\Phi^{v_{i-1}}_{\mu_{i-1}}
(C_{|M_{\mu_{i-1}}}\Phi^{v_{i}}_{\mu_{i}}-\Phi^{v_{i}}_{\mu_{i}}C_{|M_{\mu_{i}}})\Phi^{v_{i+1}}_{\mu_{i+1}}
\ldots \Phi^{v_r}_{\mu}Be^\lambda)$$ in two different ways. On one hand, we have
\begin{equation}\label{E:2}
A(\lambda,\mu)=(\Delta_{\mu_{i-1}}-\Delta_{\mu_{i}})F^{v_1,\ldots ,v_r}(\lambda,\mu).
\end{equation}
On the other hand, using the relations
\begin{equation*}
\begin{split}
(f^{(t)}_\alpha e^{(t)}_\alpha)_{|M_{\mu_{i-1}}}\Phi^{v_i}_{\mu_{i}}-
\Phi^{v_i}_{\mu_{i}}&(f^{(t)}_\alpha
e^{(t)}_\alpha)_{|M_{\mu_{i}}}\\
&=-\big\{f^{(t)}_{\alpha|M_{\mu_{i-1}}}e^{(t)}_{\alpha|V_i}+e^{(t)}_{\alpha|M_{\mu_{i-1}}}f^{(t)}_{\alpha|V_i}+(f^{(t)}_\alpha
e^{(t)}_\alpha)_{|V_i}\big\}\Phi^{v_i}_{\mu_{i}},
\end{split}
\end{equation*}
$$\rho_{|M_{\mu_{i-1}}} \Phi^{v_i}_{\mu_{i}}-
\Phi^{v_i}_{\mu_{i}}\rho_{|M_{\mu_{i}}}=-\rho_{|V_i}\Phi^{v_i}_{\mu_{i}},$$
$$x^2_{j|M_{\mu_{i-1}}}\Phi^{v_i}_{\mu_{i}}-\Phi^{v_i}_{\mu_{i}}x^2_{j|M_{\mu_{i}}}=-\big\{x^2_{j|V_i}+2x_{j|M_{\mu_{i-1}}}x_{j|V_i}\big\}\Phi^{v_i}_{\mu_{i}},$$
we deduce that
\begin{equation*}
\begin{split}
A(\lambda,\mu)=&\big\{-\sum_j
x^2_{j|V_i}-2\sum_{\alpha,t}(f_\alpha^{(t)}e_{\alpha}^{(t)})_{|V_i}-2\rho_{|V_i}
\big\}F^{v_1,\ldots ,v_r}(\lambda,\mu)\\
&+A_1(\lambda,\mu) + A_2(\lambda,\mu)
+A_3(\lambda,\mu) 
\end{split}
\end{equation*}
 where
\begin{align*}
A_1(\lambda,\mu)&=-2\sum_{j\in I} x_{j|V_i}\mathrm{Tr}(\Phi^{v_1}_{\mu_{1}} \ldots
\Phi^{v_{i-1}}_{\mu_{i-1}}x_{j|M_{\mu_{i-1}}} \Phi^{v_i}_{\mu_{i}}\ldots \Phi^{v_r}_{\mu} B e^\lambda),\\
A_2(\lambda,\mu)&=-2\sum_{\alpha,t} e^{(t)}_{\alpha|V_i} \mathrm{Tr}(\Phi^{v_1}_{\mu_1} \ldots
\Phi^{v_{i-1}}_{\mu_{i-1}}f^{(t)}_{\alpha|M_{\mu_{i-1}}} \Phi^{v_i}_{\mu_{i}}\ldots \Phi^{v_r}_{\mu} B
e^\lambda),\\ 
A_3(\lambda,\mu)&=-2\sum_{\alpha,t} f^{(t)}_{\alpha|V_i} \mathrm{Tr}(\Phi^{v_1}_{\mu_{1}} \ldots
\Phi^{v_{i-1}}_{\mu_{i-1}}e^{(t)}_{\alpha|M_{\mu_{i-1}}} \Phi^{v_i}_{\mu_{i}}\ldots \Phi^{v_r}_{\mu} B
e^\lambda).
\end{align*}
Writing $A_1$ as a sum of two equal terms:
\begin{equation*}
\begin{split}
A_1(\lambda,\mu)=&-\sum_{j\in I}
x_{j|V_i}\mathrm{Tr}(\Phi^{v_1}_{\mu_{1}}
\ldots
\Phi^{v_{i-1}}_{\mu_{i-1}}x_{j|M_{\mu_{i-1}}} \Phi^{v_i}_{\mu_{i}}\ldots \Phi^{v_r}_{\mu}
B e^\lambda)\\
&-\sum_{j\in I} x_{j|V_i}\mathrm{Tr}(\Phi^{v_1}_{\mu_{1}}
\ldots
\Phi^{v_{i-1}}_{\mu_{i-1}}x_{j|M_{\mu_{i-1}}} \Phi^{v_i}_{\mu_{i}}\ldots
\Phi^{v_r}_{\mu} B e^\lambda),
\end{split}
\end{equation*}
 and using the intertwining properties $\Phi x_j=(x_j
\otimes 1 + 1 \otimes x_j)\Phi$ and $(1 \otimes \Phi x_j)=\Phi x_j -(x_j \otimes 1)\Phi $
repeatedly in the first and second term of $A_1(\lambda,\mu)$ respectively, we get
\begin{equation}\label{E:2.5}
\begin{split}
A_1(\lambda,\mu)&=-\sum_{j\in I} x_{j|V_i}\left( \sum_{t<i}
x_{j|V_t}-\sum_{t\geq i}x_{j|V_t}\right)F^{v_1,\ldots ,v_r}(\lambda,\mu)\\ &\qquad - \sum_{j
\in I} x_{j|V_i}\left(\mathrm{Tr}(x_{j|M_{\mu_{0}}}\Phi^{v_1}_{\mu_{1}}\ldots
\Phi^{v_r}_{\mu} B e^\lambda)+\mathrm{Tr}(\Phi^{v_1}_{\mu_{1}} \ldots
\Phi^{v_r}_{\mu}x_{j|M_{\mu}} B e^\lambda)\right).
\end{split}
\end{equation}

Now 
there are two cases to consider, depending on whether $j \in I_1$ or
$j \in I_2$. In the first case,  by the cyclicity of the trace, we have
\begin{equation*}
\begin{split}
\mathrm{Tr}(x_{j|M_{\mu_{0}}}\Phi^{v_1}_{\mu_{1}}\ldots \Phi^{v_r}_{\mu} B
e^\lambda)+\mathrm{Tr}&(\Phi^{v_1}_{\mu_{1}} \ldots \Phi^{v_r}_{\mu}x_{j|M_{\mu}} B
e^\lambda)\\
=&\mathrm{Tr}(\Phi^{v_1}_{\mu_{1}}
\ldots\Phi^{v_r}_{\mu}(Bx_j+x_j)_{|M_{\mu}}Be^\lambda).
\end{split}
\end{equation*}
and
\begin{equation}\label{E:2.51}
\begin{split}
\mathrm{Tr}(\Phi^{v_1}_{\mu_{1}} \ldots
\Phi^{v_r}_{\mu}(Bx_j+x_j)_{|M_{\mu}}Be^\lambda)=&2\mathrm{Tr}(\Phi^{v_1}_{\mu_{1}}
\ldots
\Phi^{v_r}_{\mu}x_{j|M_{\mu}}Be^\lambda)\\
=&2 \frac{\partial}{\partial x_j}F^{v_1,\ldots
,v_r}(\lambda,\mu)
\end{split}
\end{equation}
where the differentiation is taken with respect to the parameter
$\lambda$. 

Let us now deal with the second case. It is easy to check using Lemma 3.2 that 
for any $x\in \l^\perp$, one has
\begin{equation*}
\begin{split}
\mathrm{Tr}(x_{|M_{\mu_{0}}}\Phi^{v_1}_{\mu_{1}}\ldots \Phi^{v_r}_{\mu} B
e^\lambda)+\mathrm{Tr}&(\Phi^{v_1}_{\mu_{1}} \ldots \Phi^{v_r}_{\mu}x_{|M_{\mu}} B
e^\lambda)\\
&=-\sum_i(C_Tx)_i\mathrm{Tr}(\Phi^{v_1}_{\mu_{1}} \ldots \Phi^{v_r}_{\mu}B
e^\lambda).
\end{split}
\end{equation*}
Therefore, we obtain
\begin{equation}\label{EE:2}
\begin{split}
A_1(\lambda,\mu)
 =\bigg(-\sum_{t<i} x_{j|V_i}x_{j|V_t}+\sum_{t\geq i}&
x_{j|V_i}x_{j|V_t}-2\sum_{j \in I_1} x_{j|V_i}\frac{\partial}{\partial
  x_j} \\
&+\sum_{j \in I_2} x_{j|V_i}\sum_{l=1}^r
C_Tx_{j|V_l}\bigg) F^{v_1,\ldots,v_r}(\lambda,\mu)\\
=\bigg(-\sum_{t<i} x_{j|V_i}x_{j|V_t}+\sum_{t\geq i}&
x_{j|V_i}x_{j|V_t}-2\sum_{j \in I_1} x_{j|V_i}\frac{\partial}{\partial
  x_j} \\
&+\sum_{j \in I_2} x_{j|V_i}\sum_{l=1,l \neq i}^r
C_Tx_{j|V_l}\bigg) F^{v_1,\ldots,v_r}(\lambda,\mu)
\end{split}
\end{equation}
where in the last equality
we used the skew-symmetry of $\sum_jC_Tx_j\otimes x_j$ to get rid of terms with 
$i=l$.

\paragraph{}We now compute $A_2(\lambda,\mu)$. Using the intertwining
property $\Phi f_\alpha=(1 \otimes f_\alpha + f_\alpha \otimes 1) \Phi$, we have
\begin{equation}\label{E:?}
\begin{split}
\mathrm{Tr}&(\Phi^{v_1}_{\mu_{1}} \ldots \Phi^{v_{i-1}}_{\mu_{i-1}}
f^{(t)}_{\alpha|M_{\mu_{i-1}}}\ldots
\Phi^{v_r}_{\mu}Be^\lambda)\\ & =\mathrm{Tr}(\Phi^{v_1}_{\mu_{1}} \ldots
f^{(t)}_{\alpha|M_{\mu_{i-2}}}\Phi^{v_{i-1}}_{\mu_{i-1}}\ldots \Phi^{v_r}_{\mu}Be^\lambda) 
+f^{(t)}_{\alpha|V_{i-1}} F^{v_1,\ldots ,v_r}(\lambda,\mu)\\ &=
\mathrm{Tr}(f^{(t)}_{\alpha|M_{\mu_{0}}}\Phi^{v_1}_{\mu_{1}} \ldots
\Phi^{v_r}_{\mu}Be^\lambda) + \big( f^{(t)}_{\alpha|V_{1}}+ \ldots +
f^{(t)}_{\alpha|V_{i-1}}\big)F^{v_1,\ldots ,v_r}(\lambda,\mu)\\ &=e^{-( \alpha,\lambda )}
\mathrm{Tr}(\Phi^{v_1}_{\mu_{1}} \ldots
\Phi^{v_r}_{\mu}Bf^{(t)}_{\alpha|M_{\mu'}}e^\lambda) + \big(
f^{(t)}_{\alpha|V_{1}}+
\ldots + f^{(t)}_{\alpha|V_{i-1}}\big)F^{v_1,\ldots ,v_r}(\lambda,\mu).
\end{split}
\end{equation}

Hence by Lemma \ref{Blem}
\begin{equation}\label{E:?1}
\begin{split}
\mathrm{Tr}(\Phi^{v_1}_{\mu_{1}}& \ldots \Phi^{v_{i-1}}_{\mu_{i-1}}
f^{(t)}_{\alpha|M_{\mu_{i-1}}}\ldots
\Phi^{v_r}_{\mu}Be^\lambda)\\
=&\bigg(\sum_{j<i}\sum_{l=0}^\infty
e^{-l(\alpha,\lambda)}B^lf_{\alpha|V_j}^{(t)}+\sum_{j\geq
i}\sum_{l=1}^\infty
e^{-l(\alpha,\lambda)}B^lf_{\alpha|V_j}^{(t)}\bigg)F^{v_1,\ldots
,v_r}(\lambda,\mu).
\end{split}
\end{equation}

A similar computation (with $e_\alpha$ moving to the right) shows that
\begin{equation}\label{E:?2}
\begin{split}
\mathrm{Tr}(\Phi^{v_1}_{\mu_{1}}& \ldots \Phi^{v_{i-1}}_{\mu_{i-1}}
e^{(t)}_{\alpha|M_{\mu_{i-1}}}\ldots
\Phi^{v_r}_{\mu}Be^\lambda)\\
=&-\bigg(\sum_{j<i}\sum_{l=1}^\infty
e^{-l(\alpha,\lambda)}B^{-l}e_{\alpha|V_j}^{(t)}+\sum_{j\geq
i}\sum_{l=0}^\infty
e^{-l(\alpha,\lambda)}B^{-l}e_{\alpha|V_j}^{(t)}\bigg)F^{v_1,\ldots
,v_r}(\lambda,\mu).
\end{split}
\end{equation}
Adding (\ref{E:?1}) and (\ref{E:?2}) and using the relation
\begin{equation}\label{E:?3}
[e_\alpha^{(t)},B^lf_\alpha^{(t)}]=\cases
(\theta^{(t)}_\alpha)^{l/N_\alpha}h_\alpha & \mathrm{if}\;
\alpha\in\langle \Gamma_3 \rangle\;\mathrm{and}\;N_\alpha|l,\\
0 & \mathrm{else}\endcases
\end{equation}
we get
\begin{equation}\label{E:I3}
\begin{split}
A_2&(\lambda,\mu)+A_3(\lambda,\mu)-2\sum_{\alpha,t}
(f_\alpha^{(t)}e_\alpha^{(t)})_{V_i} F^{v_1,\ldots
v_r}(\lambda,\mu)\\
=&-2 \bigg(\sum_{\alpha,t} \sum_{j <i}e^{(t)}_{\alpha|V_i} (T_1(\alpha))_{|V_j}
+\sum_{\alpha,t}\sum_{j
>i}e^{(t)}_{\alpha|V_i}(T_2(\alpha))_{|V_j}\bigg)F^{v_1,\ldots
,v_r}(\lambda,\mu)\\
 &+2\bigg(\sum_{\alpha,t}\sum_{j
<i}f^{(t)}_{\alpha|V_i}(T_3(\alpha))_{|V_j}F^{v_1,\ldots
,v_r}(\lambda,\mu)+\sum_{\alpha,t}\sum_{j
>i}f^{(t)}_{\alpha|V_i}(T_4(\alpha))_{|V_j}\bigg)F^{v_1,\ldots
,v_r}(\lambda,\mu)\\
&-2
\sum_{\alpha \in \langle \Gamma_3 \rangle,t}\frac{1}{1-\theta_\alpha^{(t)}e^{-N_\alpha( \alpha,\lambda )}}
h_{\alpha|V_i}F^{v_1,\ldots
,v_r}(\lambda,\mu).
\end{split}
\end{equation}
where 
$$
T_1(\alpha)=\sum_{l=0}^\infty e^{-l(\alpha,\lambda)}B^lf_\alpha^{(t)}, \qquad
T_2(\alpha)=\sum_{l=1}^\infty
e^{-l(\alpha,\lambda)}B^lf_\alpha^{(t)},$$
$$
T_3(\alpha)=\sum_{l=1}^\infty e^{-l(\alpha,\lambda)}B^{-l}e_\alpha^{(t)}, \qquad
T_4(\alpha)=\sum_{l=0}^\infty
e^{-l(\alpha,\lambda)}B^{-l}e_\alpha^{(t)}.$$

\paragraph{}Combining (\ref{E:2}), (\ref{EE:2}) and (\ref{E:I3}) now gives the following
relation for
$F^{v_1,\ldots ,v_r}(\lambda,\mu)$:
\begin{equation*}
\begin{split}
\bigg(&\sum_{j \in I_1} (x_j)_{|V_i} \otimes  \frac{\partial}{\partial x_j}+ \rho_{|V_i}-K(\lambda)_{|V_i}+\sum_{j>i}r(\lambda)_{|V_i \otimes V_j}-\sum_{j<i}r(\lambda)_{|V_j\otimes V_i}\bigg) F^{v_1,\ldots,v_r}(\lambda,\mu)\\
&=\frac{1}{2}( \Delta_{\mu-\nu_{i+1}-\ldots
-\nu_{r}}-\Delta_{\mu-\nu_i-\ldots -\nu_r})F^{v_1,\ldots ,v_r}(\lambda,\mu)
\end{split}
\end{equation*}
where 
$$K(\lambda)=
\sum_{\alpha\in\langle \Gamma_3\rangle,t}\frac{1}{1-\theta_\alpha^{(t)}e^{-N_\alpha( \alpha,\lambda )}}
h_\alpha.$$
A direct computation shows that
\begin{equation}\label{E:indiana}
\sum_{j \in I_1}x_j \frac{\partial}{\partial x_j} \mathrm{Tr}_{|M_{-\rho}}(Be^\lambda)+(\rho + K(\lambda))  \mathrm{Tr}_{|M_{-\rho}}(Be^\lambda)=0.
\end{equation}
 It is easy to deduce (\ref{E:prop1}) from the above equations.\qed

\paragraph{3.3. The second order equation for traces.}

\begin{theo}\label{SecOr} The function
$\F^{v_1,\ldots ,v_r}(\lambda,\mu)$ satisfies the following second order differential
equation :
\begin{equation}\label{E:prop12}
 \bigg(\sum_{j \in I_1} \frac{\partial^2}{\partial x_j^2}- 
\sum_{l,n=1}^r S_T(\lambda)_{|V_l \otimes V_n}\bigg)
=(\mu+\rho,\mu+\rho)\F^{v_1,\ldots ,v_r}(\lambda,\mu)
\end{equation}
where 
\begin{equation*}
\begin{split}
S_T(\lambda)=\sum_{\alpha,t}\sum_{k=0}^{\infty}
\sum_{v=1}^{\infty} e^{-(s+v)(\alpha,\lambda)} (B^sf_\alpha^{(t)}
\otimes B^{-v}e_\alpha^{(t)}+B^{-v}e_\alpha^{(t)}\otimes 
B^sf_\alpha^{(t)})\\ -\sum_{j
  \in I_2} \frac{1-C_T}{2}x_j \otimes \frac{1-C_T}{2}x_j.
\end{split}
\end{equation*}
\end{theo}
\noindent
\textit{Proof.} Consider
$$A'(\lambda,\mu)=\mathrm{Tr}(\Phi^{v_1}_{\mu_{1}} \ldots \Phi^{v_r}_{\mu}C_{|M_{\mu}}Be^\lambda).$$
On one hand, we have
$$A'(\lambda,\mu)=\Delta_{\mu} F^{v_1,\ldots ,v_r}(\lambda,\mu),$$
and on the other hand, 
\begin{equation*}
\begin{split}
A'(\lambda,\mu)=&\mathrm{Tr}\bigg(\Phi^{v_1}_{\mu_{1}} \ldots \Phi^{v_r}_{\mu} \big(\sum_{j\in I} x_j^2+ 2\rho+
2\sum_{\alpha,t} f^{(t)}_\alpha
e^{(t)}_\alpha\big)_{|M_\mu}Be^\lambda\bigg)\\
=&A'_1(\lambda,\mu)+A'_2(\lambda,\mu)+A'_3(\lambda,\mu),
\end{split}
\end{equation*}
 where
\begin{align*}
A'_1(\lambda,\mu)&=\sum_{j\in I} \mathrm{Tr}(\Phi^{v_1}_{\mu_{1}} \ldots
\Phi^{v_r}_{\mu}x^2_{j|M_\mu}Be^\lambda),\\ 
A'_2(\lambda,\mu)&=2\mathrm{Tr}(\Phi^{v_1}_{\mu_{1}} \ldots
\Phi^{v_r}_{\mu}\rho_{|M_\mu}Be^\lambda)=2\frac{\partial}{\partial \rho}F^{v_1,\ldots ,v_r}(\lambda,\mu),\\
A'_3(\lambda,\mu)&=2\sum_{\alpha,t} \mathrm{Tr}(\Phi^{v_1}_{\mu_{1}} \ldots
\Phi^{v_r}_{\mu}(f_\alpha^{(t)}e_{\alpha}^{(t)})_{|M_\mu}Be^\lambda),
\end{align*}
where the differentiation in the second equation is taken with respect to $\lambda$. To compute $A'_1(\lambda,\mu)$, note that, as in (\ref{E:2.5}), (\ref{E:2.51}),
\begin{equation}\label{E:Pip1}
A'_1(\lambda,\mu)=\sum_{j \in I_1} \frac{\partial^2}{\partial x_j^2} F^{v_1,\ldots ,v_r}(\lambda,\mu)+\sum_{j \in I_2} \mathrm{Tr}(\Phi^{v_1}_{\mu_{r-1}} \ldots \Phi^{v_r}_{\mu}x^2_{j|M_\mu}Be^\lambda).
\end{equation}
The second term on the r.h.s of (\ref{E:Pip1}) can be evaluated by the same method as in the derivation of the KZB equations: for all $j\in I_2$ we
have
\begin{equation}\label{E:Pip2}
\mathrm{Tr}(\Phi^{v_1}_{\mu_{1}} \ldots
\Phi^{v_r}_{\mu}x^2_{j|M_\mu}Be^\lambda)=\sum_{l,
s}\frac{1-C_T}{2}x_{j|V_l}\frac{1-C_T}{2}x_{j|V_s}F^{v_1,\ldots ,v_r}(\lambda,\mu).
\end{equation}
We now compute $A'_3(\lambda,\mu)$. We consider two cases\\
\textit{Case 1: $\alpha \not\in \langle \Gamma_3\rangle$.} By
the intertwining property and the cyclicity of the trace again, we have
\begin{equation}\label{E:??}
\begin{split}
\mathrm{Tr}&(\Phi^{v_1}_{\mu_{1}} \ldots
\Phi^{v_r}_{\mu}(f_\alpha^{(t)}e^{(t)}_\alpha)_{|M_\mu}Be^\lambda)\\
&=\sum_{l=1}^r
f_{\alpha|V_l}\mathrm{Tr}(\Phi^{v_1}_{\mu_{1}} \ldots \Phi^{v_r}_{\mu}e^{(t)}_{\alpha|M_\mu}Be^\lambda)+e^{-(
\alpha,\lambda)}\mathrm{Tr}(\Phi^{v_1}_{\mu_{1}} \ldots
\Phi^{v_r}_{\mu}(e_{\alpha}^{(t)}B(f^{(t)}_\alpha))_{|M_\mu}Be^\lambda).
\end{split}
\end{equation}
Applying this equation repeatedly and using the relation (\ref{E:?3}), we obtain
\begin{equation}\label{E:Pip2.5}
\mathrm{Tr}(\Phi^{v_1}_{\mu_{1}} \ldots
\Phi^{v_r}_{\mu}(f_\alpha^{(t)}e^{(t)}_\alpha)_{|M_\mu}Be^\lambda)=\sum_{l=1}^r
\sum_{s=0}^{\infty} e^{-s(
\alpha,\lambda)}B^s(f_{\alpha})_{|V_l}\mathrm{Tr}(\Phi^{v_1}_{\mu_1} \ldots
\Phi^{v_r}_{\mu}e^{(t)}_{\alpha|M_\mu}Be^\lambda).
\end{equation}
\textit{Case 2: $\alpha \in \langle \Gamma_3 \rangle$.} Applying (\ref{E:??}) $N_\alpha$ times
yields
\begin{equation*}
\begin{split}
\mathrm{Tr}(\Phi^{v_1}_{\mu_{1}} \ldots
\Phi^{v_r}_{\mu}&(f_\alpha^{(t)}e^{(t)}_\alpha)_{|M_\mu}Be^\lambda)\\
=&\sum_{l=1}^r
\sum_{s=0}^{N_\alpha-1} e^{-s(
\alpha,\lambda)}B^s(f_{\alpha})_{|V_l}\mathrm{Tr}(\Phi^{v_1}_{\mu_1} \ldots
\Phi^{v_r}_{\mu}e^{(t)}_{\alpha|M_\mu}Be^\lambda)\\
 &+e^{-N_\alpha(
\alpha,\lambda)}\theta_\alpha^{(t)}\mathrm{Tr}(\Phi^{v_1}_{\mu_{1}}
\ldots
\Phi^{v_r}_{\mu}h_{\alpha|M_\mu}Be^\lambda)\\ 
&+e^{-N_\alpha(\alpha,\lambda)}\theta_\alpha^{(t)}\mathrm{Tr}(\Phi^{v_1}_{\mu_{1}} \ldots
\Phi^{v_r}_{\mu}(f_\alpha^{(t)}e^{(t)}_\alpha)_{|M_\mu}Be^\lambda),
\end{split}
\end{equation*}
from which it follows that
\begin{equation}\label{E:Pip3}
\begin{split}
\mathrm{Tr}&(\Phi^{v_1}_{\mu_{1}} \ldots
\Phi^{v_r}_{\mu}(f_\alpha^{(t)}e^{(t)}_\alpha)_{|M_\mu}Be^\lambda)\\
=&\frac{1}{1-\theta_\alpha^{(t)}e^{-N_\alpha(
\alpha,\lambda)}}\big\{
\sum_{l=1}^r \sum_{s=0}^{N_\alpha-1} e^{-s(
\alpha,\lambda)}B^s(f_{\alpha})_{|V_l}\mathrm{Tr}(\Phi^{v_1}_{\mu_{1}}
\ldots \Phi^{v_r}_{\mu}e^{(t)}_{\alpha|M_\mu}Be^\lambda)\\
 &+\theta_\alpha^{(t)}e^{-N_\alpha(
\alpha,\lambda)}\mathrm{Tr}(\Phi^{v_1}_{\mu_{1}}
\ldots
\Phi^{v_r}_{\mu}h_{\alpha|M_\mu}Be^\lambda)\big\}.
\end{split}
\end{equation}
The formula for $A'_3(\lambda,\mu)$ now follows from  (\ref{E:?2}). Equations
(\ref{E:Pip1}),(\ref{E:Pip2}), (\ref{E:Pip3}) imply the following second-order differential
equation for $F^{v_1,\ldots ,v_r}(\lambda,\mu)$:
\begin{equation}\label{E:ohio}
\begin{split}
\bigg( \frac{1}{2}\sum_{j \in I_1} \frac{\partial^2}{\partial x_j^2} +
\frac{\partial}{\partial \rho}+H(\lambda)&\bigg)F^{v_1,\ldots
  ,v_r}(\lambda,\mu)\\
&=\frac{1}{2}\bigg(\Delta_\mu  + \sum_{l,n=1}^r S_T(\lambda)_{V_l \otimes V_n}
\bigg)F^{v_1,\ldots ,v_r}(\lambda,\mu)
\end{split}
\end{equation}
where 
$$H(\lambda)=\sum_{{\alpha}\in\langle \Gamma_3\rangle,t}\frac{1}
{1-\theta_\alpha^{(t)}e^{-N_\alpha(\alpha,\lambda)}}\frac{\partial}{\partial h_\alpha}.$$
In particular,
$$\bigg(\frac{1}{2}\sum_{j \in I_1} \frac{\partial^2}{\partial x_j^2} +
\frac{\partial}{\partial \rho}+H(\lambda)\bigg)\mathrm{Tr}_{|M_{-\rho}}(Be^\lambda)=\frac{1}{2}\Delta_\rho \mathrm{Tr}_{|M_{-\rho}}(Be^\lambda)$$
 which, together with (\ref{E:indiana}) and (\ref{E:ohio}) yields (\ref{E:prop12}). \qed

\paragraph{3.4. Diagonalization of the KZB and the second order 
operators.}

Denote by $K_j(\la)$ the differential operators appearing 
on the left hand side of the KZB equations, and by $D(\la)$ 
the second order operator appearing in Theorem \ref{SecOr}. 
These are operators on the space of functions of 
$\la$ with values in \linebreak $(V_1\otimes...\otimes V_r)^\l$. 
It is known \cite{Fel} that $K_j$ commute with each other. 
Besides, it can be shown that the operators $K_j$ commute with $D$
(in fact, this is also clear from the discussion below). 
This gives rise to the problem of simultaneous diagonalization 
of these operators. More precisely, the problem can be formulated as 
follows. 

Fix a weight $\nu\in \l^\perp$. 
Fix a generic point $\xi\in \l^*$, and 
consider the space of formal series 
$$
W_{\nu,\xi}:=
e^{(\la,\xi+\rho)}(V_1\otimes...\otimes V_r)[\nu][[e^{-(\alpha_i,\la)}]].
$$
It is clear that the operators $K_j,D$ act naturally in this space 
and are upper triangular with respect to the natural ordering. 
The problem is to find a (topological) basis of 
$W_{\nu,\xi}$ in which these operators are diagonal.

The following proposition provides such a basis. 

Let ${\mathcal B_i}$ be homogeneous bases of $V_i$, and 
${\mathcal B}(\nu)$ be the set of collections $(v_1,...,v_r)$  
of vectors $v_i\in {\mathcal B_i}$ such that the sum of their weights is 
$\nu$.

\begin{prop} For generic $\xi$, 
the functions $\F^{v_1,..,v_r}(\la,\frac{1+C_T}{2}\nu+\xi-
\sum n_i\alpha_i)$, where $n_i\ge 0$, and $(v_1,...,v_r)$ run through 
${\mathcal B}$, form a common topological eigenbasis of the operators $K_j,D$
in the space $W_{\nu,\chi}$. 
\end{prop} 

This proposition follows immediately from the theorems of this section:
the fact that the listed functions form a basis is obvious, so 
the only thing to be shown is that they are eigenfunctions, 
which was shown above. 

\paragraph{3.5. Quantum integrable systems associated to generalized 
Belavin-Drinfeld triples for simple Lie algebras.}

In the case when the Lie algebra $\g$ is finite dimensional, one can define 
other differential operators which commute with $K_j,D$. 

Namely, if $Z$ is any element of the center of $U(\g)$ then 
there exists a unique differential operator $D_Z$
on $\l^*$ with values in
 $\text{End}((V_1\otimes...\otimes V_r)[{\l^\perp}])$ such that 
$$
D_Z\F=\delta_B(\la)\mathrm{Tr}(\Phi^{v_1}_{\mu_{1}} \ldots 
\Phi^{v_r}_{\mu}Z_{|M_{\mu}}Be^\lambda).
$$
For example, $D_C=D-(\rho,\rho)$. 

It is easy to see that $D_{Z_1Z_2}=D_{Z_1}D_{Z_2}$, so 
the operators $D_Z$ form a commutative algebra. 
It is clear that these operators also diagonalize in the basis of the previous 
section. 
Thus, we get a "quantum integrable system", whose 
eigenstates are the functions $\F^{v_1,..,v_r}$. 

In the special case $\Gamma_i=\Gamma$, $T=id$, this system is 
a generalized trigonometric Calogero-Moser system considered in \cite{EFK2}. 
     
\section{Classical dynamical r-matrices with spectral parameter}
\paragraph{}Applying the construction of Section 2 to an 
(untwisted) affine Lie algebra $\hat{\g}$ and using the evaluation map $\mathrm{ev}_z:\; \hat{\g} \to \g$, one can obtain solutions of the classical dynamical Yang-Baxter equation with spectral parameter. This is done as follows.
\paragraph{} Let $\g$ be a simple complex Lie algebra, and let $\tilde{\g}=\g[t,t^{-1}] \oplus \C c \oplus \C \partial$ be the associated affine Kac-Moody algebra, where $c$ is the central element and $\partial$ is the grading element. The commutation relations in $\tilde{\g}$ are:
$$
[c,\tilde{\g}]=0,\qquad
[xt^n,yt^m]=[x,y]t^{n+m}+n\delta_{n,-m}c,\qquad
[\partial,xt^n]=nxt^n\qquad \forall\;x,y \in \g.$$
\paragraph{}Recall that the Cartan subalgebra of $\tilde{\g}$ is $\tilde{\h}=\h \oplus \C c \oplus \C \partial$. Then $\tilde{\h}^*=\h^* \oplus \C \Lambda_0 \oplus \C \delta$ where $\Lambda_0,\delta$ are defined by $\langle\Lambda_0,\h\rangle =\langle \Lambda_0,\partial \rangle =\langle \delta,\h \rangle=\langle \delta,c\rangle =0$ and $\langle \delta,\partial\rangle=\langle \Lambda_0,c\rangle=1$. Under the standard bilinear form on $\tilde{\h}$, $c$ and $\partial$ are orthogonal to $\h$ and we have $( c,\partial )=1$, $( c,c ) = ( \partial,\partial) =0$.\\
\hbox to1em{\hfill}The root system of $\tilde{\g}$ is $\tilde{\Delta}=(\Delta + \Z\delta) \cup
\Z^*\delta$. The root subspace corresponding to $\alpha +k\delta$ is spanned by $e_\alpha t^k$ if
$\alpha \in \Delta^+$, $f_{-\alpha}t^k$ if $\alpha \in -\Delta^+$ and 
equals $\h t^k$ if $\alpha=0$.
The system of positive roots is $\tilde{\Delta}^+=\Delta^+ \cup \N \delta \cup (\Delta + \N
\delta)$.
\paragraph{}We will consider a twisted version of the affine Lie algebra $\tilde{\g}$. Let $g$ be the dual
Coxeter number of $\g$ and set $\epsilon=e^{\frac{2i\pi}{g}}$. Consider the
automorphism $\gamma=Ad(e^{2i\pi \rho /g})$ of $\g$. We have $\gamma
(e_\alpha)=\epsilon^{|\alpha|}e_\alpha, \gamma (f_\alpha)=\epsilon^{-|\alpha|}f_\alpha,$
and $\gamma_{|\h}=Id$. Let $\tilde{\g}_\gamma$ be the subalgebra of $\tilde{\g}$
consisting of all elements $a(t)+\lambda c + \mu \partial$ satisfying $a(\epsilon
t)=\gamma(a(t))$. The elements $e_\alpha t^{|\alpha|+mg},\;f_\alpha
t^{-|\alpha|+mg},\;x_it^{mg}, c, \partial$ for $\alpha \in \Delta^+,m \in \Z$ and $(x_i)$ an
orthonormal basis of $\h$, form a $\C$-basis of
$\tilde{\g}_\gamma$. The proof of the following lemma is straightforward.
\begin{lem}\label{L:1}
The map $\phi: \tilde{\g}_\gamma \to \tilde{\g}$ defined by
$$\phi(e_\alpha t^{|\alpha|+mg})=e_\alpha t^m,\qquad \phi(f_\alpha t^{-|\alpha|+mg})=f_\alpha t^m,\qquad\phi(\partial)=g\partial +\rho$$
$$\phi(x_it^{mg})=x_i t^m,\;\;(m \neq 0),\qquad
\phi(x_i)=x_i-(\rho,x_i)\frac{c}{g},\qquad \phi(c)=\frac{c}{g}$$ is a Lie algebra isomorphism.\end{lem}

Let $(\Gamma_1,\Gamma_2,\tau)$ be a generalized Belavin-Drinfeld triple for 
$\hat\g_\gamma$. Let $\tilde \l$ be the subalgebra of $\tilde \h$ 
of elements $x$ such that $(\alpha,x)=(T(\alpha),x)$ for 
$\alpha\in \Gamma_1$. Let $\l=\tilde l\cap \h$. It is clear that 
$\tilde\l$ contains $c$ and $\partial$, so $\tilde\l=\C c\oplus 
\C\partial\oplus \l$. 

\paragraph{}For any $z \in \C^*$ let $\mathrm{ev}_z:\; \g[t,t^{-1}] \to \g$ be the evaluation map defined by $\mathrm{ev}_z(xt^n)=z^n x$ for all $x \in \g$, $n \in \Z$. 

Fix a complex number $\tau$ with positive imaginary part. 
Let $\tilde\la=\la+2\pi i\tau\delta/g$. Define 
$$
\overline{r}_T(\la,z)=(\text{ev}_z\otimes \text{ev}_1)
(\phi^{-1}\otimes \phi^{-1})(r_T(\tilde\lambda)),
$$
where $z\in \C^*$, $\lambda\in \l^*$. 

{\bf Remark.} Although the evaluation maps are not defined 
on $\partial$, this definition makes sense, since $\partial$ 
occurs in $r_T$ in a combination $c\otimes\partial +\partial \otimes c$, 
and $ev_z(c)=0$ for any $z$. 

It is clear that $\overline{r}_T$ is a Laurent series in $z$ whose 
coefficients are meromorphic functions on $\l^*$. 

\begin{prop} The series 
$\overline{r}_T(\lambda,z): \l^*\to (\g \otimes \g)^\l$ is 
convergent in a nonempty annulus, and extends to a meromorphic function 
on $\l^*\times \C^*$. Moreover, this meromorphic function satisfies 
the dynamical Yang-Baxter equation with spectral parameter:
\begin{equation*}
\begin{split}
\sum_i& x_i^{(1)} \frac{\partial}{\partial x_i}r^{23}(\lambda,\frac{z_2}{z_3})-x_i^{(2)} \frac{\partial}{\partial x_i}r^{13}(\lambda,\frac{z_1}{z_3})+x_i^{(3)}\frac{\partial}{\partial x_i}r^{12}(\lambda,\frac{z_1}{z_2})\\
&+[r^{12}(\lambda,\frac{z_1}{z_2}),r^{13}(\lambda,\frac{z_1}{z_3})]+[r^{13}(\lambda,\frac{z_1}{z_3}),r^{23}(\lambda,\frac{z_2}{z_3})]+[r^{12}(\lambda,\frac{z_1}{z_2}),r^{23}(\lambda,\frac{z_2}{z_3})]=0,
\end{split}
\end{equation*}
where $x_i$ is an orthonormal basis of $\h$.
\end{prop}

This proposition follows easily from Proposition 2.1.

Now we would like to compute $\overline{r}_T(\la,z)$ explicitly. 
For the sake of simplicity we will restrict ourselves to
triples of the form $(\Gamma_1=\Gamma_2=\Gamma,T)$ (i.e $T$ is
an automorphism of the Dynkin diagram): these triples give rise to
elliptic dynamical r-matrices. 

{\bf Remark.} More general triples give rise to
partially trigonometric and partially elliptic r-matrices. 

Let $T$ be an  automorphism of the Dynkin diagram $\Gamma$ of
$\tilde{\g}$ of order $N$. As in Section 2, let $B$ be the lift of $T$ to an automorphism of
$\tilde{\g}$ of order $N$. By Lemma~\ref{L:1}, $B$ defines an automorphism of
$\tilde{\g}_\gamma$. Note that the action of
$B$ on
$\tilde{\g}_\gamma$ preserves the principal gradation, i.e the exists a unique automorphism $\beta$ of $\g$ such
that
$B(x t^m)=\beta(x) t^m$ for any $x t^m \in \tilde{\g}_\gamma$. Furthermore, $c$ and
$\partial$ are $B$-invariant, and $\h$ is $B$-stable. 
Like before, we will choose an orthonormal basis
$x_i$ of $\h$ which is compatible with $\l=\h^B$ and $\l^\perp$. 

Let
$$\theta(u,\tau)=-\sum_{j=-\infty}^{\infty} e^{\pi i (j +\frac{1}{2})^2\tau
+ 2\pi i(j+\frac{1}{2})(u+\frac{1}{2})}$$
be the Jacobi theta function. For brevity we will not write the 
dependence on $\tau$ explicitly. Introduce the functions
$$\sigma_w(u)=\frac{\theta(w-u)\theta'(0)}
{\theta(w)\theta(u)},\qquad
\chi(u)=\frac{\theta'(u)}{\theta(u)},$$
where $\theta'$ is the derivative of $\theta$.

Let $z=e^{2\pi iu/g}$.

\begin{prop} The function $\overline{r}_T(\lambda,z): \l^*\times \C^*
\to (\g \otimes \g)^\l$ is equal to
\begin{equation*}
\begin{split}
\overline{r}_T&(\lambda,z)=\\
&-\sum_{\alpha >0} \sum_{l=0}^{N-1} \frac{1}{2\pi i}
e^{-l(\alpha,\tilde\la)+2\pi i|\alpha|u/g}
\sigma_{\frac{N}{2\pi i}(\alpha,\tilde\la)}
(u-l\tau|N\tau) e_\alpha \otimes \beta^l(f_\alpha)\\
&-\sum_{\alpha >0} \sum_{l=0}^{N-1} \frac{1}{2\pi i}
e^{l(\alpha,\tilde\la)-2\pi i|\alpha|u/g}
\sigma_{-\frac{N}{2\pi i}(\alpha,\tilde\la)}
(u-l\tau|N\tau) f_\alpha \otimes \beta^l(e_\alpha)\\
&-\sum_{i \in I}\sum_{l=0}^{N-1} \left(
\frac{1}{2}\delta_{l0}+\frac{1}{2\pi i}\chi(u-l\tau|N\tau) \right)x_i \otimes
\beta^l x_i +\frac{1}{2}\sum_{i \in I_2}
\frac{\beta+1}{\beta-1}x_i \otimes x_i
\end{split}
\end{equation*}
(here $\alpha$ runs through positive roots of $\g$). 
\end{prop}

The proof of this proposition is by a direct calculation. 

\paragraph{}Thus we see 
that the above construction gives a dynamical r-matrix with the number of dynamical
parameters equal to $d-1$, where $d$ is the number of orbits of $T$ on the Dynkin diagram.
\paragraph{Examples.} i) When $T=1$, this yields 
(up to a gauge transformation) Felder's elliptic r-matrix, as shown in \cite{EV1}, Section 4.6.\\
ii) Let $\g=sl(n)$ and let $T$ be the rotation of the Dynkin diagram by $2\pi
k/n$. Note that in this case (and in this case only) $d=1$ and we obtain a
non-dynamical r-matrix with spectral parameter. It is easy to check that it is equal to Belavin's classical elliptic r-matrix with modulus $\tau$
(see \cite{BD}).

\section{Twisted traces of intertwiners for affine Lie algebras}
\paragraph{}In this section, we apply the same procedure as in Section 3 in the case of an affine Lie algebra $\tilde{\g}$, but we consider finite-dimensional (evaluation) modules $V_i$ rather than modules from the category
$\mathcal{O}$. 
We keep the notations of
Section 4.

\paragraph{}The Lie algebra $\tilde{\g}_\gamma$ has a
triangular decomposition
$\tilde{\g}_\gamma=\tilde{\n}^+_\gamma \oplus \tilde{\h}
\oplus \tilde{\n}^-_\gamma$ where $\tilde{\n}^+_\gamma$
(resp. $\tilde{\n}^-_\gamma$) is the subalgebra spanned by
$e_\alpha t^{|\alpha|+mg}, \;m \geq 0$, $f_\alpha
t^{-|\alpha|+mg},\;x_it^{mg},$ $ \;m>0$ (resp.
spanned by $f_\alpha t^{-|\alpha|+mg},\;m \leq 0$,
$e_\alpha t^{|\alpha|+mg},\;x_it^{mg},
\;m<0$). This decomposition allows us to define
highest weight Verma modules and dual Verma modules in
the usual way. 
We put $\tilde{\mu}=\mu + \frac{k}{g}\Lambda_0 -\frac{g(\mu,\mu)}{2(k+g)}\delta$ and denote by $M_{\mu,k}$ the Verma module $M_{\tilde{\mu}}$. 

\paragraph{}Finally, we define evaluation representations of $\tilde{\g}_\gamma$. Let $V$ be a finite-dimensional highest weight $\g$-module with highest weight vector $v_0$. For any $\Delta \in \C$ let $z^{-\Delta}V[z,z^{-1}]$ denote the $\tilde{\g}$-evaluation module,
with $\tilde{\g}$-action given by
$$xt^n.(vP(z))=xvz^nP(z),\qquad c.(vP(z))=0,\qquad \partial.(vP(z))=vz\frac{d}{dz}P(z)$$
for all $x \in \g,\; v \in V,\; P(z) \in z^{-\Delta}\C[z,z^{-1}].$ Let $z^{-\Delta}V_\gamma[z,z^{-1}] \subset z^{-\Delta}V[z,z^{-1}]$ be the subspace of all (Laurent) polynomials $v(z)$ satisfying $v(\epsilon z)=e^{2i \pi \rho /g}v(z)$. Then $z^{-\Delta}V_\gamma[z,z^{-1}]$ is a $\tilde{\g}_\gamma$-module (which is isomorphic to the usual $\tilde{\g}$-evaluation module). From now on, we simply write $z^{-\Delta}V[z^{\pm 1}]$ for $z^{-\Delta}V_\gamma[z^{\pm 1}]$.

\paragraph{}Proposition 3.1 admits an analogue in this situation (see e.g \cite{EFK1},\cite{Et}) :
\begin{prop} Let $\mu \in \h^*$ and $k \in \C$ and let $V$ be a finite-dimensional $\g$-module. Suppose that $M^*_{\mu,k}$ is irreducible. Then for each $v \in V[\lambda-\mu]$ there exists a unique $\tilde{\g}_\gamma$ intertwiner 
$$\tilde{\Phi}^v_{\mu,k}(z): M_{\lambda,k} \to M_{\mu,k} \hat{\otimes} z^{-\Delta}V[z^\pm 1]$$
where
 $\Delta=g\frac{(\lambda,\lambda)-(\mu,\mu)}{2(k+g)}$ such that
$$\langle v_{\mu,k}^*, \tilde{\Phi}^v_{\mu,k}(z)v_{\lambda,k}\rangle= z^{-\Delta}v.$$ 
\end{prop}
Here $\hat\otimes$ denotes the completed tensor product.
\paragraph{}Let $V_1,\ldots V_r$ be finite-dimensional $\g$-modules, and let $v_1\in V_1,\ldots v_r \in V_r$ be homogeneous vectors of weight $\nu_1,\ldots \nu_r$ respectively. Set $\nu=\sum \nu_i$, 
and assume $\nu\in\l^\perp$. We will consider composition of intertwining operators :
$$
\tilde{\Phi}^{v_1}_{\mu-\nu_2-\ldots -\nu_{r},k}(z_1) \ldots \tilde{\Phi}^{v_r}_{\mu,k}(z_r) B e^\lambda :
\;M_{\mu',k}
\to  M_{\mu-\nu,k} \otimes z_1^{-\Delta_1}V_1[z_1^{\pm 1}] \otimes \ldots \otimes z_r^{-\Delta_r}V_r[z_r^{\pm
1}]$$ where $\lambda \in (\h^*)^B,\;\mu \in \h^*$, where 
$(\tilde\mu',\alpha)=(\tilde\mu,T\alpha)$ for $\alpha\in\Gamma_1$, and
$\Delta_i=g\frac{2(\nu_{i},\mu-\nu_{i+1}-\ldots -\nu_{r})-(\nu_i,\nu_i)}{2(k+g)}$. This composition lives in a certain
completion (see \cite{EFK1}, $\S 3$) of the space
$$\mathrm{Hom}\big(M_{\mu',k},M_{\mu-\nu,k}\otimes V_1 \otimes \ldots \otimes V_r\big)\otimes 
z_1^{-\Delta_1}\ldots z_r^{-\Delta_r}\C[z_1^{\pm 1},\ldots,z_r^{\pm 1}].$$ Let us fix some
$\tau \in \C$. For $\lambda \in \h^*$ set
$\tilde{\lambda}=\la+2\pi i
\tau
\delta/g$.  Consider the trace function
$$F^{v_1,\ldots v_r}_k(\lambda,\mu,\mathbf{z})=\mathrm{Tr}\;(\tilde{\Phi}^{v_1}_{\mu-\nu_2-\ldots
-\nu_{r},k}(z_1)
\ldots
\tilde{\Phi}^{v_r}_{\mu,k}(z_r) B e^{\tilde{\lambda}})$$
where $\mathbf{z}=(z_1,\ldots ,z_r)$.
 Finally, let
$$\hat{\delta}_B(\lambda)=\big(\mathrm{Tr}_{|M_{0,k}}(Be^{\tilde{\lambda}})\big)^{-1}
$$
and set
$$\F^{v_1,\ldots v_r}_k(\lambda,\mu,\mathbf{z})=\hat{\delta}_B(\lambda)F_k(\lambda,\mu,\mathbf{z}).$$
Our main results in the affine case are 
\begin{theo} The function $\F^{v_1,\ldots v_r}_k(\lambda,\mu,\mathbf{z})$ satisfies the following system of
differential equations for $i=1,\ldots r$ :
\begin{equation}\label{E:theoaff}
\begin{split}
-\frac{k + g}{g}z_i &\frac{\partial}{\partial z_i}\F^{v_1,\ldots v_r}_k(\lambda,\mu,\mathbf{z})= \sum_{j \in I_1}x_{j|V_i}
\frac{\partial}{\partial x_j}
\F^{v_1,\ldots v_r}_k(\lambda,\mu,\mathbf{z})\\
&+\big(\sum_{j>i}\overline{r}_T(\lambda,z_i/z_j)_{|V_i \otimes V_j} -
\sum_{j<i}\overline{r}_T(\lambda,z_j/z_i)_{|V_j \otimes V_i}
\big)\F^{v_1,\ldots v_r}_k(\lambda,\mu,\mathbf{z}).
\end{split}
\end{equation}
\end{theo}

These equations are called the KZB equations (see \cite{Fel}). 

\begin{theo}
The function $\F^{v_1,\ldots ,v_r}(\lambda,\mu,\mathbf{z})$ satisfies
the following second order differential equation :
\begin{equation}\label{E:prop22}
-\frac{k+g}{\pi i }\frac{\partial}{\partial
\tau}\F_k(\lambda,\mu,\mathbf{z})
=\sum_{j \in I_1}\frac{\partial^2}{\partial
x_j^2}\F^{v_1,\ldots v_r}_k(\lambda,\mu,\mathbf{z})
-\sum_{i,j=1}^r\overline{S}^{ij}(\lambda,z_i/z_j)
\F^{v_1,\ldots v_r}_k(\lambda,\mu,\mathbf{z}),
\end{equation}
where 
$$\overline{S}_T(\la,z):=(\text{ev}_z\otimes \text{ev}_1)(\phi^{-1}
\otimes \phi^{-1})(S_T(\la)), 
$$
and $S_T(\la)$ is the function defined in Theorem 3.2 (for $\tilde \g$). 
\end{theo}

The proofs of these theorems are completely parallel to the proofs of the 
theorems of Section 3. 

{\bf Remark 1.} If $\Gamma_1=\Gamma_2=\Gamma$ 
then $\overline{S}_T(\lambda,z)$ can be expressed in terms of elliptic 
functions. 

{\bf Remark 2.} In section 3 we saw that the differential operator 
\linebreak $D=\Delta_\l-S_T(\la)$ (where $\Delta_l$ is the Laplacian of $\l$)
on functions with values in the space of weight zero under $\l$ 
in some representation of the Lie algebra can be included in a quantum 
integrable system. A similar statement holds for the 
operator $D_\tau=\Delta_\l-\overline{S}_T(\la,1)$, 
whose coefficients are elliptic (if $\Gamma_i=\Gamma$). 
To obtain other operators commuting with $D_\tau$, 
it is necessary to apply the construction of Section 3.5 
 to central elements of a completion of $U(\tilde \g)$
at the critical level $k=-g$, as explained in \cite{EFK2}.   
In the case $T=1$, this gives a generalization of the elliptic
Calogero-Moser system (see \cite{EFK2}).

{\bf Remark 3.} Although we considered
 only untwisted affine algebras, the results can be easily
generalized to the twisted case. 

\paragraph{Acknowledgments.} P.E. would like to thank  IHES for hospitality.
O.S. is grateful to the Harvard Mathematics Dept. for the kind invitation. The work of P.E was
partially supported by the NSF grant 9700477, and was partly done when P.E. was employed by the Clay Mathematical Institute as a CMI prize fellow.
P.E. is grateful to Travis Schedler for useful discussions.

\end{document}